\documentclass[12pt]{amsart}
\usepackage{amsfonts,amsmath,amsthm,amscd,amssymb,latexsym,cite,verbatim,texdraw,floatflt,
caption2}
\RequirePackage{hyperref}

\usepackage[a4paper]{geometry}

\theoremstyle
{plain}

\begin{document}

\title{Closeness and linkness in balleans}

\author{ I. Protasov, K. Protasova}

\maketitle
\vskip 5pt

{\bf Abstract.} A set $X$ endowed with a coarse structure is called ballean or coarse space. 
For a ballean  $(X, \mathcal{E})$,  we say that two subsets  $A$, $B$  of $X$   are close (linked) if there exists an entourage $E\in \mathcal{E}$ 
 such that $A\subseteq E [B]$,  $B\subseteq  E[A]$  (either $A, B$ are bounded or contain unbounded close subsets). 
 We explore the following general question:  
which information about a ballean is contained and can be extracted  from the relations of closeness and linkness. 

\vspace{6 mm}


\vspace{3 mm}


\vspace{6 mm}

\section{Introduction}

Given a set $X$, a family $\mathcal{E}$ of subsets of $X\times X$ is called a {\it coarse structure }  on $X$  if
\vskip 7pt

\begin{itemize}
\item{}   each $E\in \mathcal{E}$  contains the diagonal  $\bigtriangleup _{X}$,
$\bigtriangleup _{X}= \{(x,x)\in X: x\in X\}$;
\vskip 5pt

\item{}  if  $E$, $E^{\prime} \in \mathcal{E}$ then $E\circ E^{\prime}\in\mathcal{E}$ and
$E^{-1}\in \mathcal{E}$,   where    $E\circ E^{\prime}=\{(x,y): \exists z((x,z) \in  E,  \   \ (z, y)\in E^{\prime})\}$,   $E^{-1}=\{(y,x): (x,y)\in E\}$;
\vskip 5pt

\item{} if $E\in\mathcal{E}$ and $\bigtriangleup_{X}\subseteq E^{\prime}\subseteq E  $   then
$E^{\prime}\in \mathcal{E}$;
\vskip 5pt

\item{}   $\bigcup  \mathcal{E}= X\times X $.

\end{itemize}
\vskip 7pt

A subfamily $\mathcal{E}^{\prime} \subseteq \mathcal{E}$  is called a
{\it base} for $\mathcal{E}$  if,
 for every $E\in \mathcal{E}$, there exists
  $E^{\prime}\in \mathcal{E}^{\prime}$  such  that
  $E\subseteq E ^{\prime}$.
For $x\in X$,  $A\subseteq  X$  and
$E\in \mathcal{E}$, we denote
$E[x]= \{y\in X: (x,y) \in E\}$,
 $E [A] = \bigcup_{a\in A}   \      E[a]$
 and say that  $E[x]$
  and $E[A]$
   are {\it balls of radius $E$
   around} $x$  and $A$.

The pair $(X,\mathcal{E})$ is called a {\it coarse space} \cite{b13}  or a ballean \cite{b10}, \cite{b12}.

For a ballean $(X,\mathcal{E})$, a  subset $B \subseteq X$   is called {\it bounded} if $B \subseteq  E[x]$  for some
$E\in \mathcal{E}$ and $x\in X$.
A ballean $(X,\mathcal{E})$ is called {\it bounded (unbounded)} if $X$ is bounded (unbounded).
The family $\mathcal{B}_{(X, \mathcal{E})}$
 of all bounded subsets of  $(X, \mathcal{E})$
is called the {\it bornology} of $(X, \mathcal{E})$.
We denote by  $\mathcal{B}^{\sharp}_{(X,\mathcal{E})}$  the set of all ultrafilters  $\varphi$  on $X$  such  that $X\setminus B \in \varphi$  for each $B\in \mathcal{B}_{(X,\mathcal{E})}$.
\vskip 10pt

{\bf Definition 1. } Given a  ballean  $(X,\mathcal{E})$,  we say that  subsets  $A$, $B$  of $X$   are 

\vskip 7pt

\begin{itemize}
\item{}  {\it close} (write $A\delta B$) if there exists $E \in \mathcal{E}$ such  that $A\subseteq  E [B]$,  $B\subseteq  E[A]$; 

\item{}  {\it  linked }  (write $A\lambda   B$) if either $A, B$ are  bounded or there exist unbounded  subsets $A^{\prime} \subseteq  A$,  $B^{\prime} \subseteq  B$ such that  $A^{\prime}\delta   B^{\prime}$.

\end{itemize}

The relations  $\delta$  and  $\lambda $,  called  {\it closeness}  and  {\it linkness}, are special  kinds of  asymptotic  proximities  [8]. 
We  note that  $\delta$  and  $\lambda $
   play the  central  part in ultrafilter description of the  Higson's  corona 
 $\nu  (X,  \mathcal{E})$ of   $(X,  \mathcal{E})$,  see [7] 
 and Section 6.
  The  negation of  $\lambda$ (namely,  asymptotic disjointness) is used for definition of normal  balleans, see [6]  and Section 4.  
The relation $\delta$   is an equivalence and each   $\delta$-class is a   connected  component  in the hyperballean  of $(X, \mathcal{E})$, 
 see  [4], [11].

In above and other cases,  $\delta$  and  $\lambda $ were used as a tool for studing balleans,  but
$\delta$  and  $\lambda $  are interesting   for their own sake.  
In this paper, we concentrate around the following general  questions:  
 \vskip 7pt

{\it Which properties of a ballean  $(X, \mathcal{E})$ can be recognized by 
$\delta _{(X,  \mathcal{E})}$  or  $\lambda _{(X,  \mathcal{E})}$ ?   
How can one detect if a ballean $(X, \mathcal{E})$  is uniquely  determined by 
  $\lambda _{(X,  \mathcal{E})}$
 or $\delta _{(X,  \mathcal{E})}$? }
\vskip 7pt

To make these questions more precise, we need some definitions. 
\vskip 10pt

{\bf Definition 2. } We  say that two balleans $(X, \mathcal{E})$, 
$(X, \mathcal{E^{\prime}})$
 are
 $\delta$-{\it equivalent  ($\lambda$-equivalent)} if 
$\delta _{(X,  \mathcal{E})}=$
  $\delta _{(X,  \mathcal{E^{\prime}})}$
($\lambda _{(X,  \mathcal{E})}=$
  $\lambda _{(X,  \mathcal{E^{\prime}})}$).
\vskip 7pt

Every $\delta$-equivalent balleans  are  $\lambda$-equivalent, but not vice versa, see Section   2.

If   $(X, \mathcal{E})$, 
$(X, \mathcal{E^{\prime}})$,   are $\lambda$-equivalent then
 $ \mathcal{B}_{(X,\mathcal{E})}= \mathcal{B}_{(X,\mathcal{E^{\prime}})}$. 
\vskip 10pt

{\bf Definition 3. } A  ballean  $(X,\mathcal{E})$ is called {\it  $\delta${-rigid}  ({$\lambda$-rigid})}  if, 
 for every coarse structure  
$\mathcal{E^{\prime}}$ on $X$ such that
$\delta _{(X,  \mathcal{E})}= \delta _{(X,  \mathcal{E^\prime})}$
($\lambda _{(X,  \mathcal{E})}= \lambda _{(X,  \mathcal{E^\prime})}$),
  we have 
$ \mathcal{E}  = \mathcal{E^\prime} $ . 
\vskip 7pt

Thus, every  $\delta${-rigid}  ({$\lambda$-rigid})    ballean  $(X,\mathcal{E})$ is uniquely determined by the relation
  $\delta _{(X,  \mathcal{E})} \ ( \lambda _{(X,  \mathcal{E})})$.
Every $\lambda$-rigid  ballean  is  $\delta$-rigid.

For  balleans  $(X,  \mathcal{E})$, uniquely   determined by 
$ \mathcal{B_{ (X, E)}}$ see [3].
\vskip 10pt

{\bf Definition 4. } We  say that a class $ \mathcal{K}$ of balleans  is 

\vskip 7pt

\begin{itemize}

\item{}  {\it $\delta$-stable  ($\lambda$-stable), for any two   $\delta$-equivalent ($\lambda$-equivalent) balleans 
 $(X,\mathcal{E})$,   $(X,\mathcal{E}^{\prime})$,
 if    $(X,\mathcal{E})\in \mathcal{K}$ then 
 $(X,\mathcal{E}^{\prime})\in \mathcal{K}$;}

\item{}  {\it $\delta$-rigid  ($\lambda$-rigid), for any two   $\delta$-equivalent ($\lambda$-equivalent) balleans 
 $(X,\mathcal{E})$,   $(X,\mathcal{E}^{\prime})$
   $\in  \mathcal{K}$, we have
 $\mathcal{E} = \mathcal{E}^{\prime}$.}

\end{itemize}

By [8, Theorem 3.4],  for every ballean  $(X,\mathcal{E})$, there exists the strongest by inclusion coarse structure $\mathcal{E}_{\delta} $ 
 on $X$  such that  
 $(X,\mathcal{E})$,  $(X,\mathcal{E}_{\delta})$
 are  $\delta$-equivalent. 
Analogously, there exists the strongest by inclusion coarse structure $\mathcal{E}_\lambda$ on  $X$  such that
$(X,\mathcal{E})$,  $(X,\mathcal{E}_{\lambda})$
 are  $\lambda$-equivalent.
\vskip 10pt

{\bf Definition 5. } We  say that a   ballean  $(X,\mathcal{E})$  is $\delta$-{\it strong} ($\lambda$-{\it strong})  if  $ \mathcal{E}= \mathcal{E}_{\delta}$  ($ \mathcal{E}= \mathcal{E}_{\lambda}))$. 
\vskip 10pt

In Section 2 
we prove that every discrete ballean is $\lambda$-rigid, every coarsely  discrete   ballean is 
$\delta$-rigid but needs not to be  $\lambda$--rigid.
\vskip 10pt

In Section 3
 we show that the class of metrizable  balleans is $\lambda$-rigid but not  $\lambda$-stable, the class
 of submetrizable balleans is $\lambda$-stable.
 Every metrizable ballean  $(X, \mathcal{E})$ is $\lambda$-strong (so $\delta$-strong)  but $(X, \mathcal{E})$ is $\lambda$-rigid if and only if $(X,\mathcal{E})$ is discrete. 
 Is every metrizable ballean $\delta$-rigid? This question remains open. 
\vskip 10pt

In Section 4 we prove that the class of normal  balleans is  $\lambda$-stable,  but not    $\delta$-rigid. 
\vskip 10pt

In Section 5 we prove that the class of finitary group  balleans is  $\lambda$-rigid.
\vskip 10pt

In Section 6 
we show that the parallelity relation on 
$ \mathcal{B}_{ (X, \mathcal{E})} ^{\sharp}$ 
uniquely defines  $\delta_{ (X, \mathcal{E})}$ and apply this statement to   balleans $(X,\mathcal{E})$
 with finite  
$ \mathcal{B}_{ (X, \mathcal{E})} ^{\sharp}$. 
We conclude the paper with an ultrafilter characterization  of the normality. 

%2 zzzzzzzzzzzzzzzzzzzzzzzzzzzzzzzzzzzzzzzzzzzzzzzzzzzzzzzzzzzzz

\section{Discrete balleans}

We  recall that a ballean $(X, \mathcal{E})$ is {\it discrete} if, for every $E\in \mathcal{E}$, 
there exists a bounded subset $B$ such that $E[x]=\{ x\}$ 
 for each $x\in X \backslash B$. A subset $Y\subseteq X$  is called {\it discrete} if the subballean
 $(Y,  \mathcal{E} |_{ Y})$ is  discrete. 

\vskip 10pt

{\bf Theorem 1.} {\it For a ballean $(X, \mathcal{E})$,  the following statements are equivalent: 
\vskip 10pt

$(1)$   $(X, \mathcal{E})$  is discrete; 
\vskip 7pt

$(2)$	if  $Y\subseteq X$, $Z\subseteq  X$ and $Y\delta Z$ then $Y\setminus Z$ and   $Z\setminus Y$ are bounded;    
\vskip 7pt

$(3)$	 if  $Y\subseteq X$, $Z\subseteq  X$ and $Y\lambda Z$ then either  $Y,   Z$ are  bounded or  
$Y\cap Z$ is unbounded;    
\vskip 10pt

Proof. } The implications  $(1)\Longrightarrow   (2)\Longrightarrow  (3)$ are evident. 

We assume that (3) holds but $(X,  \mathcal{E})$ is not discrete.
 Then there exists $E\in  \mathcal{E}$
 such that the set  
  $\{ x\in X: |E[x]|>1    \}$ is unbounded. 
 We choose an unbounded subset $Y$ of $X$ such that 
$|E[y]|>1$ for each $y\in Y$ and $E[y]\cap  E[y^{\prime}]= \emptyset$  
for all distinct 
$y, y^{\prime} \in Y$.
 For each $y\in Y$, we pick 
$z_y \in E[y] \setminus  \{ y\}$ and denote
 $Z=\{ z_y: y\in Y\}$. Clearly, $Y\lambda Z$ but $Y\cap Z =\emptyset  $ contradicting (3).    $ \ \  \  \Box$ 
\vskip 10pt

{\bf Corollary 2.}  {\it Every discrete ballean is $\lambda$-rigid.}  

\vskip 10pt

To continue, we need the following construction from [5]. 

Let $(X, \mathcal{E})$ be an  unbounded ballean and $\varphi$ be a filter on $X$ such that  $X\setminus B\in \varphi $ for each  $B\in  \mathcal{B}_{ (X, \mathcal{E})}  $.
We denote  by $ \mathcal{E}_{\varphi}$ a coarse structure on $X$ with the base 
 $\{ H_{(E,\phi)}: E\in  \mathcal{E}$, $\Phi\in \varphi \}$, where
\vskip 10pt

%$$$$$$$$$$$$$$$$$$formula

\[  H_{(E,\phi)} [x]  = \left\{
\begin{array}{ll}
\{x\} & \textrm{if  } x\in \Phi\textrm{,}\\
E[x]\setminus \Phi  & \textrm{if }  x\in X\setminus \Phi\textrm{.}
\end{array} \right. \]

%%%%%%%%%%%%%%%%%

\vskip 10pt

Clearly,   $ \mathcal{E}_{\varphi}\subseteq \mathcal{E}$ and $\mathcal{B}_{ (X, \mathcal{E})}  = \mathcal{B}_{ (X,  \mathcal{E}_{\varphi})}$.

\vskip 10pt

{\bf Theorem 3.} {\it Let $(X, \mathcal{E})$ be  a ballean,  $E_{0}\in \mathcal{E}$.  Assume that there exists an unbounded subset $Y$  of $X$ such  that    $|E_{0}[y]| >1$,  $y\in Y$  and every unbounded subset  of $Y$ can be partitioned in two unbounded subsets. 
Then there exists an ultrafilter $\varphi$  on $X$  such that 
$ \mathcal{E}_{\varphi}\subset \mathcal{E}$ and 
$ (X, \mathcal{E}),$   $  (X,  \mathcal{E}_{\varphi})$ are  $\lambda$-equivalent, so $ (X, \mathcal{E})$  is not  $\lambda$-rigid.

\vskip 10pt

Proof.}  We  take an arbitrary  ultrafilter $\varphi \in \mathcal{B}^{\sharp}_{ (X, \mathcal{E})} $  such that  $Y\in \varphi$.
By the assumption,  $ \mathcal{E}_{\varphi}\subset \mathcal{E}$.

We take arbitrary subsets $A, B$ of $X$  linked in  $(X, \mathcal{E})$ and show that $A, B$  are linked in $(X, \mathcal{E}_{\varphi})$. 
If $A, B$  are bounded in  $(X, \mathcal{E})$  then  $A, B$  are bounded in  $(X, \mathcal{E}_{\varphi})$,  so we  suppose that 
$A, B$  are unbounded  in  $(X, \mathcal{E})$.
We choose unbounded subset $A^{\prime}\subseteq A$,  subset $B^{\prime}\subseteq B$ close in   $ (X, \mathcal{E})$ . If $A^{\prime} \cap Y $ and $B^{\prime} \cap Y $ are bounded then $A^{\prime},B^{\prime}$ are close in $(X, \mathcal{E}_{\varphi})$ by the construction of  $\mathcal{E}_{\varphi}$.  We suppose that $A^{\prime} \cap Y$ is unbounded and partition $A^{\prime} \cap Y$ in two unbounded subsets $C, D$.  Since $\varphi$ is an ultrafilter, either  $C\notin \varphi$ or $D\notin \varphi$. Hence, either $C$  or $D$ is linked to $B^{\prime}$ in $(X, \mathcal{E}_{\varphi})$.  
$ \ \  \  \Box$ 
\vskip 10pt

Let  $(X, \mathcal{E})$ be a ballean.  We recall that a subset $Y$ of $X$  is {\it large}
  if there exists $E\in \mathcal{E}$ such that $X= E[Y].$

Clearly, $Y$  is large if  and only if  $Y\delta X$. We note that every  unbounded subset of $X$ is linked with $X$.

A ballean   $(X, \mathcal{E})$ is called {\it coarsely discrete} if  $X$  contains a large discrete subset.

% Since every $\lambda$-rigid  ballean is $\delta$-rigid, by corollary 2,  the class of coarsely rigid  balleans is $\delta$-stable. 

\vskip 10pt

%%%%%%%%%%%%%%%%%%%%%%%%%%%%%

{\bf Theorem 4.} {\it Every coarsely discrete ballean $(X, \mathcal{E})$  is $\delta$-rigid.
\vskip 10pt

Proof. }  We suppose that  $\mathcal{E}=\mathcal{E}_{\delta} $, take  a  $\delta$-equivalent ballean  $(X, \mathcal{E}^{\prime})$ 
and show that  $\mathcal{E}^{\prime}= \mathcal{E}$.

Let $Y$  be a large discrete subset of  $(X, \mathcal{E})$,  $E\in  \mathcal{E}$ 
and $E[Y] =X$. We take an arbitrary $E^{\prime}\in  \mathcal{E}^{\prime}$ such  that $E^{\prime}[Y] =X$ and  prove that there exists a bounded subset $B$ of $X$  such that  $E^{\prime}[x] =E[x]$ for each $x\in X\setminus B$.

We choose $H\in  \mathcal{E}$  such that $E\subseteq H$,  $E^{\prime}\circ   E^{\prime} \subseteq H$ and 
take a bounded subset $C$  such that $H[y]\cap H[y^{\prime}] = \emptyset $ for all distinct $y, y^{\prime} \in Y\setminus C$. 
Then   $E^{\prime} [y]\subseteq E[y] $  and  
 $E^{\prime} [z]\subseteq E[y] $
for all  $y\in Y\setminus C$, $z\in E[y]$. 
Since  $E^{\prime} [Y]=X$,  there exists a bounded subset $B$  such that  $C\subseteq B$  and  
$E^{\prime} [x]= E[x] $ for each  $x\in X\setminus B$.  $ \ \  \  \Box$ 
\vskip 10pt

\vskip 10pt

{\bf Example 5.}
We define a coarse structure  $ \mathcal{E}$ on the  set $ \mathbb{N}$ of natural numbers  such that the ballean    $(\mathbb{N}, \mathcal{E})$  is  $\delta$-rigid  but   not  $\lambda$-rigid.  We  denote  $A=\{ (2n+1, 2n+2),   (2n+2, 2n+1):  n\in \omega \}$ 
and put  $\mathcal{E} =  \{( F \times F) \cup A:  F\in   \mathbb{N}_{<\omega} \}.$
Since the set  $2 \mathbb{N}$ is  large in  $(\mathbb{N}, \mathcal{E})$,  by Theorem 4, $(\mathbb{N}, \mathcal{E})$  is $\delta$-rigid.  By Theorem 3 with $Y=2 \mathbb{N}$,  $(\mathbb{N}, \mathcal{E})$ is not  $\lambda$-rigid.

%%%%%%%%%%%%%%%%%%%%%%%%%%%%

\section{Ordinal and submetrizable balleans}

We recall that a ballean $(X, \mathcal{E})$  is {\it ordinal} if $\mathcal{E}$  has a base linearly ordered  by inclusion. In this case  $\mathcal{E}$  has a base well-ordered by inclusion. 

A ballean $X, \mathcal{E}$  is called {\it metrizable} if there exists a metric $d$ on $X$  such that  $\mathcal{E}$   has  the  base  $\{ \{ (x, y): d(x,y)<r \} : r\in \mathbb{R}^{+} \}$.
A ballean $(X, \mathcal{E})$ is metrizable if and only if $\mathcal{E}$  has a countable base [12, Theorem 2.1.1].  Hence, every metrizable ballean is ordinal. 

We denote by $\mathcal{M}$   and $\mathcal{L}$    the cases of metrizable  and ordinal balleans. 

\vskip 10pt

{\bf Theorem 6.} {\it The following statements hold: 

\vskip 10pt

$(1)$  every ordinal ballean is $\lambda$-strong; 
\vskip 7pt

$(2)$	 $\mathcal{M}$   and $\mathcal{L}$   are  $\lambda$-rigid;    
\vskip 7pt

$(3)$	  an    ordinal ballean    $(X, \mathcal{E})$          is $\lambda$-rigid if and only if $(X, \mathcal{E})$  is discrete;
\vskip 10pt

$(4)$	 $\mathcal{M}$   and $\mathcal{L}$   are not   $\lambda$-stable;    
\vskip 7pt

Proof. }
$(1)$   let   $(X, \mathcal{E})$   be an unbounded ordinal ballean, $\{ L_{\lambda}: \lambda< \kappa \} $ be a base of  $ \mathcal{E}$
well-ordered by inclusion.      We take a ballean   $(X, \mathcal{E}^{\prime})$     $\lambda$-equivalent to   $(X, \mathcal{E})$  and show that $ \mathcal{E}^{\prime} \subseteq  \mathcal{E} $.   Assume the contrary and choose $ E^{\prime} \in \mathcal{E}^{\prime} \setminus  \mathcal{E} $. Then we can choose a $\kappa$-sequence  $(x_{\alpha})_ {\alpha<\kappa}$ in $X$  such that 
\vskip 7pt

\begin{itemize}
\item{} $ L _{\alpha}[x_{\alpha}] \cap  L _{\beta}[x_{\beta}]=  \emptyset$,  $\alpha < \beta< \kappa$;

\item{} $E^{\prime}[x_{\alpha}] \cap  E^{\prime}[x_{\beta}]=  \emptyset$,  $\alpha < \beta< \kappa$;

\item{} $ E^{\prime}[x_{\alpha}] \setminus  L _{\alpha}[x_{\alpha}] \neq    \emptyset$,  $\alpha < \kappa$.

\end{itemize}

For each  $\alpha < \kappa$, we pick  
$y _{\alpha } \in  E^{\prime}  [x_{\alpha}] \setminus  L _{\alpha}[x_{\alpha}]$ and put 
 $X^{\prime}  = \{ x_{\alpha}  : \alpha < \kappa\}$, 
 $Y  = \{y_{\alpha}  : \alpha < \kappa\}$.
Then $X^{\prime} $,  $Y$ are close in  $(X, \mathcal{E} ^{\prime} )$ . For each $E\in  \mathcal{E}$,  there 
exists $\beta<\kappa$  such that 
$ E[x_{\alpha}] \subseteq  L _{\alpha}[x_{\alpha}]$,   $\alpha < \beta$. 
Hence, 
 $E [X^{\prime}] \cap  Y$  is bounded and $X^{\prime}$, $Y$ are  not linked in  $\mathcal{E}$ contradicting $\lambda$-equivalence of   $(X, \mathcal{E})$  and   $(X, \mathcal{E}^{\prime})$ .

\vskip 7pt
$(2)$   follows from $(1)$. 
\vskip 7pt

$(3)$   We note that every unbounded subset of an ordinal ballean can be partitioned in two unbounded subsets and apply Theorem 3.
\vskip 7pt

$(4)$  follows from  Theorem 3.  
$ \ \  \  \Box$ 
\vskip 10pt

By Theorem 4, the statement  $(3)$  of  Theorem 6 does not hold for $\delta$ in place of $\lambda$.

\vskip 10pt

{\bf Question 1.} {\it Is every  metrizable ballean  $\delta$-rigid?   Equivalently, is  $\mathcal{M}$   $\delta$-stable?}

\vspace{3 mm}

We  recall that an unbounded ballean  $(X, \mathcal{E})$  is {\it submetrizable}  if there exists a coarse structure  $\mathcal{E}^{\prime}$   on $X$  such that  $\mathcal{E}  \subseteq   \mathcal{E}^{\prime}  $
and   $(X, \mathcal{E}^{\prime})$  is unbounded and metrizable. 

For a ballean $(X, \mathcal{E})$,  a  function   $f: X  \longrightarrow \mathbb{R}$  
is called {\it macro-uniform}  if,  for every  
 $E\in \mathcal{E},  $ there exists  $r_E \in \mathbb{R}^{+}$   such that 
$diam \  f (E[x]) \leq r _E$
for each $x\in X$. 
We denote $mu (X, \mathcal{E} )$ the  family of all   macro-uniform  functions on $X$. 

By [12, Theorem  2.2.3],  an unbounded ballean   $(X, \mathcal{E})$ is submetrizable  if and only if there exists an unbounded function  $f\in mu  (X, \mathcal{E})$.

\vskip 10pt
%%%%%%%%%%%%%%%%%%%%%%%%%%

{\bf Theorem 7.} {\it If balleans   $(X, \mathcal{E})$, $(X, \mathcal{E^{\prime}})$ are $\lambda$-equivalent  then 
 $ mu(X, \mathcal{E})=  mu(X, \mathcal{E^{\prime}})$. 

\vskip 7pt

Proof. }
We assume that there exists  $f\in    mu (X, \mathcal{E}^{\prime}) \setminus   mu(X, \mathcal{E})$  and choose 
 $E\in \mathcal{E}$ and a sequence  $(y_n)_{n\in \omega}$ in $X$ such that  $diam \  f(E[y_n])> 2n$. For each $n\in \omega$,  we pick  $z_n \in  \  E[y_n]$ such  that $|f (y_n)|- f (z_n)|> n$.
Passing to subsequences,  we may  suppose that 
 $|f (y_n)|- f (z_m)|> n$  for all $m\geq n$.  We denote $Y= \{ y_n : n\in \omega \} $,   $Z= \{z_n : n\in \omega \} $ and note that 
$Y \delta _{(X, \mathcal{E})} Z$.
Since $f\in mu \ (X,  \mathcal{E}^{\prime})$, for every $  E^{\prime}\in \mathcal{E}^{\prime}$, there exists $k\in \mathbb{N}$
such that  $diam \  f(\mathcal{E}^{\prime} [y_n] )< k$
for all $n\in \omega$. 
Then $\mathcal{E}^{\prime} [Y] \cap Z$ is finite, so $Y, Z$ are not linked in  $(X,  \mathcal{E}^{\prime})$.
$ \ \  \  \Box$ 
\vskip 10pt

{\bf Corollary 8.} {\it The class of submetrizable balleans is   $\lambda$-stable. }

\vspace{3 mm}

A ballean  $(X, \mathcal{E})$ is called {\it mu-bounded} if  $(X, \mathcal{E})$ is not submetrizable. For  mu-bounded  balleans, see [1].
By   Corollary 8, the class of  
mu-bounded  balleans is  $\lambda$-stable. 

%%%%%%%%%%%%%%%%%%%%%%%%%%%%%%%%%%%%

\section{Normal balleans}

For a  ballean  $(X, \mathcal{E})$, two subsets $A, B$  of $X$ are called {\it asymptotically disjoint} if, for every $E\in \mathcal{E}$, 
$E[A] \cap E[B]$  is bounded. We note that   unbounded subsets $A, B$ are asymptotically disjoint  if and only if $A, B$ are not linked. 

A subset $U$  of $X$ is called an {\it asymptotic  neighbourhood} 
of a subset $A$  if , for every $E\in \mathcal{E}$,  $E[A] \setminus U$ is bounded. 

A ballean  $(X, \mathcal{E})$ is called {\it normal}  [6]  if any two asymptotically disjoint subsets have disjoint asymptotic neighbourhoods.

A function $f: X\longrightarrow  \mathbb{R}$ is called {\it slowly oscillating} if, for any $E\in \mathcal{E}$  and $\varepsilon >0$,  there exists a bounded subset $B$  of $X$  such that $diam \ f(E[x])<\varepsilon $ for each $x\in X$.
By [6, Theorem 2.2], a ballean $(X, \mathcal{E})$ is normal if  and only if, for any two disjoint and asymptotically  disjoint subsets $A, B$  of $X$, there exists a slowly oscillating function $f: X\longrightarrow  [0,1]$ such that $f|_ A =0$,  $f|_B =1$.

We denote by $so (X, \mathcal{E})$ and $sob (X, \mathcal{E})$ the families of  all  and all bounded slowly oscillating  functions on $X$.

\vskip 10pt
%%%%%%%%%%%%%%%%%%%%%%%%%%

{\bf Theorem 9.} {\it If balleans   $(X, \mathcal{E})$, $(X, \mathcal{E^{\prime}})$ are $\lambda$-equivalent  then 
 $ sob \  (X, \mathcal{E})= sob \  (X, \mathcal{E^{\prime}})$. 

\vskip 7pt

Proof. } We suppose  the contrary  and let $f\in \  sob \  (X,\mathcal{E^{\prime}})\setminus sob \  (X,\mathcal{E})$. 
We denote $\mathcal{B} = \mathcal{B} _{(X, \mathcal{E})} = \mathcal{B} _{(X, \mathcal{E^{\prime}})}$
and take $\varepsilon >0$ and $E\in \mathcal{E}$   such that,  for every $B\in \mathcal{B}$, there exists $x_{B} \in X\setminus B$  such that $diam \ f(E [x_{B}])> \varepsilon$.  We pick $y_{B}\in E [x_{B}]$ such that 
 $|f(x_{B}) - f(y_{B})|>\varepsilon $.
Then we enumerate  $\mathcal{B} = \{B_{\alpha} : \alpha < \kappa  \}$, denote $Z= \{x_B : B\in  \mathcal{B} \}$
and define  a mapping $h: Z\longrightarrow  X$  as follows. 

We put  $h(x_{B_0})=y_{B_0}$ and let $h$  is already defined for each $x_{B_{\alpha}}$, $\alpha< \beta$. If $h$  is not defined at 
  $x_{B_{\beta}}$ then we put
  $f(x_{B_{\beta}}) = y _{B_{\beta}}$. 

We take an ultrafilter $p\in \mathcal{B} ^{\sharp }$ such that $Z\in p$, denote $q=h^{\beta}(p)$ and observe that 
 $|f^{\beta}(p)- f^{\beta}(q)|\geq\varepsilon $.
We choose $P\in p$, $Q\in p$  such that $P\subseteq  Z$,  $Q\subseteq  Z$,   
$|f (x)- f^{\beta}(p)|  < (\varepsilon/2) $, $x\in P$  and  $|f (x)- f^{\beta}q)|  < (\varepsilon/2) $, $x\in Q$. 
By the construction  of $h$, we have $P\delta  _{(X, \mathcal{E})}   Q$, so 
$P\lambda  _{(X, \mathcal{E})}   Q$. 
By the assumption,  $P\lambda  _{(X, \mathcal{E}^{\prime})}   Q$. 
Hence. there exist $ E^{\prime}\in  \mathcal{E}^{\prime}$ and an unbounded $P^{\prime} \subseteq P$ such that $E^{\prime}[x]\cap Q \neq\emptyset $ for each $x\in P^{\prime}$.  Since $f$  is slowly oscillating  in  $(X, \mathcal{E^{\prime}})$, there  exists 
$B\in  \mathcal{B}$  such that  $diam \ f(E^{\prime} [x])< (\varepsilon/2)$ for each $x\in X\setminus B$.  Then 
$ f(E^{\prime} [x])\cap Q= \emptyset$  for each $x\in P^{\prime}  \setminus B$ and we get a contradiction with the choice of $P^{\prime}$.
$ \ \  \  \Box$ 
\vskip 10pt

{\bf Corollary 10.} {\it The class of normal balleans is $\lambda$-stable. }
\vskip 10pt

If  $(X, \mathcal{E})$, $(X, \mathcal{E^{\prime}})$  are  metrizable and 
 $sob(X, \mathcal{E}) = sob(X, \mathcal{E^{\prime}})$ then 
 $\lambda_{(X, \mathcal{E})}=\lambda_{(X, \mathcal{E^{\prime}})}$ 
and,  by Theorem 6(2),
  $\mathcal{E}=  \mathcal{E^{\prime}}$.

\vskip 5pt

We show that the class of normal balleans is not $\delta$-rigid.
\vskip 10pt

{\bf Example 11.}  We denote by $S_{\omega}$ the group of all permutations of $\omega$, by $\mathcal{E} $ the coarse structure on $\omega$ with the base
$$
\{\{ (x, y): y\in Fx \}  \cup \triangle_\omega  : F\in  [S _{\omega}]^{<{\omega}}\}
$$
and take  the strongest coarse structure $\mathcal{E}^\prime $  on  $\omega$  such that
  $\mathcal{B}_{(\omega, \mathcal{E} ^{\prime})  } = [\omega]^{<\omega}$,
see Example 3 in [3]. 
Then $\mathcal{E}\subset  \mathcal{E} ^{\prime}$ and every  infinite subset of $\omega$ is large in 
 $(\omega, \mathcal{E} )$  and 
 $(\omega, \mathcal{E}^{\prime} )$. 
Hence,     
$\delta_{(\omega, \mathcal{E})}=\delta_{(\omega, \mathcal{E}^{\prime} )}$
and 
$(\omega, \mathcal{E} )$, 
 $(\omega, \mathcal{E}^{\prime} )$ are normal because any two infinite subsets of $X$ are not asymptotically  disjoint in 
 $(\omega, \mathcal{E} )$  and 
 $(\omega, \mathcal{E}^{\prime} )$.
\vskip 10pt

{\bf Question 2.}  {\it  Is $so(X, \mathcal{E}) = so(X, \mathcal{E^{\prime}})$ for any  $\lambda$-equivalent ballean 
$(X, \mathcal{E})$, $(X, \mathcal{E^{\prime}})$?}

\vskip 10pt

Following [2], we say that a ballean $(X, \mathcal{E})$  has 
{\it  bounded growth}
 if,  there exists a mapping  $f: X\longrightarrow  \mathcal{B} _{(X, \mathcal{E})} $ such that $x\in f(x)$ and, for every  
$E \in  \mathcal{E}$,  the set  $X\setminus \{ x\in X: E[x] \subseteq f(x) \} $ is  bounded. 
\vskip 10pt

{\bf Question 3.} {\it Is the class of balleans of bounded grouth $\lambda$-stable?}

%%%%%%%%%%%%%%%%%%%%%%%%%%%%

\section{Finitary balleans}

For a group $G$,  we denote by $\mathcal{F}_G $ the  coarse structure on $G$  with the base
$$
\{\{ (x, y)\in  G \times   G : y\in Fx \}  \cup \triangle_G : F\in  [G]^{< \omega}\}
$$
and  say that  $(G, \mathcal{F}_G)$ is the {\it finitary group ballean } of  $G$. The following theorem show that the class of finitary group balleans is  $\lambda$-rigid.

\vskip 10pt
%%%%%%%%%%%%%%%%%%%%%%%%%%

{\bf Theorem 12.} {\it  Let $G, H$ be group on a set $X$ such that  
 $\lambda_{(X, \mathcal{F}_G) }   =  \lambda _{(X, \mathcal{F}_H)}$.  
Then $ \mathcal{F}_G = \mathcal{F}_H$. 
\vskip 7pt

Proof. }  If $X$  is countable then the  statement  follows  from Theorem 6(2)
 because  $ \mathcal{F}_G, \   \mathcal{F}_H$ have countable bases, so we suppose that $X$ is uncountable. 

We denote  by $ \  \cdot   , \   \ast $ the group operations in  $G$,  $H$  and, on the contrary, assume that 
$ \mathcal{F}_G \setminus  \mathcal{F}_H \neq \emptyset $.   Then there  exists  $F\in  [X]^{<\omega}$ such that,  for any  $H\in  [X]^{<\omega}$,  there exists  $x\in X$   such  that  $F  \cdot    x\setminus H   \ast   x  \neq \emptyset $.

We use the following  observation: 

\vskip 7pt

$(\star)$   for every  countable  subset $Y$  of $X$,  there exists  a countable subset $Z$  such that  $Y \subset Z$  and,  for every  $A\in [Y] ^{<\omega} $ ,  there exists $z\in Z$  such that  $F \cdot  z \setminus A  \ast   z \neq  \emptyset $.

We fix an arbitrary countable subset $Y _0 $  of $X$ containing $F$ and apply $(\star)$  to get a countable subgroup  $Z_0$  of $G$,
$Y_0 \subset Z_0$. 
Then we denote by $Y_1$   the  subgroup of  $H$   generated by  $Z_0$.   Analogously,  for    $Y_1$  we use  $(\star)$   to get the subgroup  $Z_1$ 
of $G$,     $Y_1  \subset  Z_1$ , and denote by $Y_2$  the subgroup of $H$  generated by $ Z_1$.  After  $\omega$  steps, we put  $S= \bigcup _{n\in\omega}  Y_n$  and note that   $S$  is a subgroup  of $G$   and  $S$  is a  subgroup of $H$.

By the construction of $S$,  for every  $A\in [S]^{<\omega}$,  there  exists $x\in S$   such that   $F \cdot   x\setminus A\ast  x \neq \emptyset $.  Hence,  the restriction of  $\mathcal{F}_G$
 and  $\mathcal{F}_H$  to $S$  are distinct.  Since $S$  is  countable, we get a  contradiction to Theorem 6(2).  $ \ \  \  \Box$ 
\vskip 10pt

A ballean $(X, \mathcal{E})$ is called  {\it  finitary} if,  for every  $E\in \mathcal{E}$,  there exists $n\in  \omega$  such that  $|E[x]|< n$   for each  $x\in X$.  By [9],  for every finitary ballean   $(X, \mathcal{E})$,  there exists  a group  $S$  of permutations of  $X$  such  that  $ \mathcal{E} $  has  the base 
$$
\{\{ (x, y)  :   x\in Fy \}  \cup \triangle_X : F\in  [S]^{< \omega}\}.
$$

Example 11 shows that  the  class of  finitary   balleans is  not  $\delta$-stable.

 \vskip 10pt

{\bf Question 4.} {\it  Is the class of  finitary  balleans 
$\delta$-rigid? }
 \vskip 10pt

The product of two $\lambda$-rigid   balleans  needs not  to be $\lambda$-rigid: take two  metrizable  unbounded discrete balleans and apply Corollary 2 and Theorem 3.

 \vskip 10pt

{\bf Question 5.} {\it  
Is the product of two $\delta$-rigid   balleans  $\delta$-rigid ? }
 \vskip 10pt

A ballean  $(X, \mathcal{E})$ is called   {\it cellular} if $ \mathcal{E}$ has a base consisting of equivalence relations. 
Equivalently [3, Theorem 3.1.3],   $(X, \mathcal{E})$ i  is cellular if   $asdim   (X, \mathcal{E}) = 0$.

 \vskip 10pt

{\bf Question 6.} {\it  
Is the class of cellular balleans $\lambda$-stable?}

%%%%%%%%%%%%%%%%%%%%%%%%%%%%

\section{Ultrafilters}

Let  $(X, \mathcal{E})$  be a ballean. We  endow $X$  with the discrete topology and consider  the  $Stone-\check{C}ech \  \ compactification$  $\beta X $  of $X$.  We  take the points of $\beta X$ to be the ultrafilters on $X$.  Then  $\mathcal{B}^\sharp  _{(X, \mathcal{E})}$   is  a closed  subset of  $\beta X$.

Given any 
 $r, q \in  \mathcal{B}^{\sharp} _{(X, \mathcal{E})}$,  we say   that  $r, q$  are {\it parallel} 
 (and write   $r|| _{(X, \mathcal{E})} q $) if there  exists  $E\in \mathcal{E} $  such that $E[P]\in  q$   for each  $P\in p$. 
By [ 6,  Lemma 4.1], 
  $|| _{(X, \mathcal{E})}  $ is an equivalence on 
 $  \mathcal{B}^{\sharp} _{(X, \mathcal{E})}$.    We denote by  $\sim  _{(X, \mathcal{E})}  $  the minimal (by  inclusion)  closed  (in   $  \mathcal{B}^{\sharp} _{(X, \mathcal{E})}  \times     \mathcal{B}^{\sharp} _{(X, \mathcal{E})}$) equivalence  such that 
$ || _{(X, \mathcal{E})} \  \subseteq \   \sim  _{(X, \mathcal{E})}  $.
The quotient  $\nu (X, \mathcal{E})$ is called the  Higson's corona of  $(X, \mathcal{E})$. 

By  [7,  Proposition 1],   $r \sim  _{(X, \mathcal{E})} q  $ if and only if  $h^{\beta}(p)= p^{\beta}(q)$  for  every  slowly  oscillating  function $h: X\longrightarrow   [0,1]$.

\vskip 10pt
%%%%%%%%%%%%%%%%%%%%%%%%%%

{\bf Theorem 13.} {\it  Let $ \mathcal{E},  \mathcal{E}^{\prime}$ be  coarse structures on a set $X$  such that  
$  \mathcal{B}^{\sharp} _{(X, \mathcal{E})}  =     \mathcal{B}^{\sharp} _{(X, \mathcal{E}^{\prime})}$
and 
$ || _{(X, \mathcal{E})}  =    ||_{(X, \mathcal{E}^{\prime})}$.
Then  
$ {(X, \mathcal{E})},     {(X, \mathcal{E}^{\prime})}$ are $\delta$-equivalent.  If 
$  \mathcal{B}^{\sharp} _{(X, \mathcal{E})}  $
is finite then $\mathcal{E}=\mathcal{E}^{\prime}$.
\vskip 7pt

Proof. }  We suppose that there exist   $P\subseteq X$,  $Q \subseteq X$ such  that  $P \delta _{(X, \mathcal{E})} Q $
and  $P\setminus  E^{\prime} (Q)\neq \emptyset   $  for each  $E^{\prime}\in  \mathcal{E}^{\prime}$. 

Since  $P \delta _{(X, \mathcal{E})} Q $,   there exists   $E\in  \mathcal{E}$  such that  $P\subseteq E[Q]$  and  $E=E^{-1}$.
For each  $x\in P$,  we pick  $f(x)\in Q$  such that  $(x, f(x))\in E. $  If 
$p\in   \mathcal{B}^{\sharp} _{(X, \mathcal{E})}$  and  $P\in p$ then 
$ p  || _{(X, \mathcal{E})}  f^ {\beta}(p) $   and  $Q\in f^ {\beta}(p) $. 
Thus, each  ultrafilter 
$p\in   \mathcal{B}^{\sharp} _{(X, \mathcal{E})}$   such that $P\in p$  is parallel in 
$  \mathcal{B}^{\sharp} _{(X, \mathcal{E})}$ to some ultrafilter  $q$ such that  $Q\in q$. 

We take an arbitrary  ultrafilter  $r$  such that $P\setminus  E^{\prime}[Q]\in  r$ for  each  $ E^{\prime} \in  \mathcal{E}. $
Then  $P\in r$  and  $r$  is not parallel in  $\mathcal{B}^{\sharp} _{(X, \mathcal{E}^{\prime})}$ to each ultrafilter $q$   such that $Q\in q$. Hence,  $ || _{(X, \mathcal{E})} \neq  || _{(X, \mathcal{E}^{\prime})}$.

To prove the second  statement,  it  suffices  to
show that  $(X, \mathcal{E})$ is coarsely discrete and 
  apply Theorem 4.

We choose a representative from  each class of parallel ultrafilter in $  \mathcal{B}^{\sharp} _{(X, \mathcal{E})}  $.
Let $p_1,  \dots , p_m$ be the set of  obtained  representatives. 
Let   $E\in  \mathcal{E}$. 
We choose $P_1 \in p_1, \dots, P_n \in p_n$ such that  $E[P_i] \cap  P_j  = \emptyset$.
Then 
$P_1 \cup, \dots, P_n $  is  discrete and large  in  $ (X, \mathcal{E})$,  so  $ (X, \mathcal{E})$ is coarsely discrete.  
$ \ \  \  \Box$ 
\vskip 10pt

%%%%%%%%%%%%%%%%%%%%%%%%%%

{\bf Remark  14.}  Let   $\mathcal{B}^{\sharp} _{(X, \mathcal{E})}=$  $\mathcal{B}^{\sharp} _{(X, \mathcal{E}^{\prime})}$.
Does    $\delta _{(X, \mathcal{E})}=$  $\delta _{(X, \mathcal{E}^{\prime})}$
implies  
 $|| _{(X, \mathcal{E})}=|| _{(X, \mathcal{E}^{\prime})}$? The  answer  to this  question  is  negative.  
We take the balleans  
 $(\omega,  \mathcal{E}), (\omega, \mathcal{E}^{\prime})$ from Example 11 and  note that   $p,q\in \omega^{\ast} $ are parallel  in   $\mathcal{B}^{\sharp} _{(\omega,  \mathcal{E})}$ if and only if $p=f ^{\beta}(q)$ for some  bijection
 $f: \omega\longrightarrow  \omega$. 
On the other hand, we partition  $\omega =  \bigcup_{n\in \omega} W_n$ so that $|W_n|= n+1$ and take a mapping 
$h: \omega\longrightarrow  \omega$ such that $h(W_n) =  \{ a_n \}$, 
$  a_n  \in  W_n$.  We choose $p\in \omega^{\ast}$  such  that  $p$  and 
 $h^{\beta}(p)$  
 are not  isomorphic. 
By the choice of $h,  p$ and $h^{\beta} (p)$ are parallel in 
  $\mathcal{B}^{\sharp} _{(X, \mathcal{E}^{\prime})}$
but not in 
$\mathcal{B}^{\sharp} _{(X, \mathcal{E})}$.

\vskip 10pt
%%%%%%%%%%%%%%%%%%%%%%%%%%

{\bf Theorem 15.} {\it  If   $ (X,  \mathcal{E})$,   $ X,  \mathcal{E}^{\prime})$  are finitary balleans  and
$ ||_{(X,  \mathcal{E}) }=  || _{(X,  \mathcal{E}^{\prime} )} $  then
  $ \mathcal{E}=  \mathcal{E}^{\prime}$.
\vskip 7pt

Proof. }  
Let  $G, G^{\prime} $   be the  groups  of  permutations of $X$ defining  $  \mathcal{E},  \mathcal{E}^{\prime}$. 
We assume that  $  \mathcal{E}\setminus   \mathcal{E}^{\prime} \neq \emptyset $. 
Then there  is
$H\in [G]^{<\omega}$
 such that,  for every 
$A\in [G^{\prime}]^{<\omega}$,
  the set 
$Y_A  = \{ x\in X: Hx \setminus   Ax \neq \emptyset \}$
  is infinite. 
We take 
$ p\in \mathcal{B}^{\sharp}_{   (X,  \mathcal{E})}$ 
 such that
$Y_A \in p$
  for each 
$A\in [G^{\prime}]^{<\omega}$.
Let 
$H= \{h_1 , \dots , h_n \}$.
We take
$g_1 , \dots , g_n \in  G^{\prime}$ 
 such that
$h_1 p= g_1 p, \dots , h_n p = g_n p$
  and put  
$B=\{ g_1 , \dots , g_n \}$.
Then
$\{ x\in X: h_1 x =  g_1 x , \dots , h_n x =  g_n x \}\in p$
  and
$Y_B \notin p $,  contradict the choice of $ p$.
$ \ \  \  \Box$ 
\vskip 10pt

In light of Theorem 15, to answer  Question 4  affirmatively, it suffices to show that 
$\delta_{(X,\mathcal{E} )} = \delta_{(X,\mathcal{E}^{\prime} )} $
 implies
$ ||_{(X,  \mathcal{E}) }=  || _{(X,  \mathcal{E}^{\prime} )} $ for any finitary  balleans
  $ \mathcal{E},  \mathcal{E}^{\prime}$.

\vskip 10pt
%%%%%%%%%%%%%%%%%%%%%%%%%%

{\bf Theorem 16.} {\it  
For a normal  ballean  $(X,\mathcal{E} )$   and
$p, q\in  \mathcal{B}_{ (X, \mathcal{E})} ^{\sharp}$, 
 the following  statements are equivalent: 

\vskip 7pt
$(1)$ $p\sim q$;

\vskip 7pt
$(2)$  $P\lambda Q$  for each  $P\in p$,  $Q\in q$;

\vskip 7pt
$(3)$ $(p, q)\in cl ||$.

\vskip 7pt

Proof. }   $(1)\Longrightarrow (2)$. 
If  $P, Q$  are unbounded and asymptotically  disjoint then,  by the normality of
 $(X,\mathcal{E} )$, 
   there  exists a slowly  oscillating function
 $f: X\longrightarrow [0,1] $  such that  
$f|_P = 1$, $f|_Q = 0$,
so
$f^{\beta} (p)\neq  f^{\beta} (q)$
 and $p,q$ are not equivalent.
\vskip 7pt

 $(2)\Longrightarrow (3)$. 
We  choose  
$E\in \mathcal{E} $, 
unbounded  subsets 
 $P^{\prime} \subseteq P$,  $Q^{\prime} \subseteq Q$   such that there is a bijection  
$h: P^{\prime} \longrightarrow Q^{\prime}$
satisfying 
$(x, h(x))\in E$, $x\in P^{\prime}$.
If  
$ r\in \mathcal{B}_{ (X, \mathcal{E})} ^{\sharp}$ 
and
$P^{\prime} \in r$
 then
$r || h ^{\beta} (r)$
 and 
$P\in r$, $Q\in  h^{\beta} (r)$.

The implication 
 $(3)\Longrightarrow (1)$
follows from definitions of
$\sim$
  and  $||$.
$ \ \  \  \Box$

\vskip 10pt
%%%%%%%%%%%%%%%%%%%%%%%%%%

{\bf Theorem 17.} {\it  
A ballean $(X,\mathcal{E} )$    is normal if and only if $\sim = cl ||$.

\vskip 7pt

Proof. }  
By Theorem 16, if $(X,\mathcal{E} )$    is normal then  $\sim = cl ||$.
We assume that
$\sim = cl ||$
  and show that any two unbounded asymptotically disjoint subsets $P, Q$ of $X$ have disjoint  asymptotic neighborhoods.

Let $p, q\in \mathcal{B}_{ (X, \mathcal{E})} ^{\sharp}$,  $P\in p$,  $Q\in q$. 
Since  $P, Q$  are    asymptotically disjoint, we have  $(p, q)\notin cl ||$, so there exists a slowly oscillating function 
$f: X\longrightarrow   [0,1]$  such that $f^{\beta}(p)=1$, $f^{\beta}(q)=0.$ 
We denote 
$P^{\prime}= \{ x\in P: f(x)> 3/4 \}$,  $Q^{\prime}= \{ x\in Q: f(x)< 1/2 \}$. 
Then $P^{\prime} \in p$,  $Q^{\prime}\in q$  and $P^{\prime}, Q^{\prime} $ have disjoint asymptotic neighborhoods.

We fix  $p\in  \mathcal{B}_{ (X, \mathcal{E})} ^{\sharp}$ 
  such that $P\in p $  and apply above paragraph  to find  $q_1, \dots , q_n\in  \mathcal{B}_{ (X, \mathcal{E})} ^{\sharp}$,  $Q_1\in q_1, \dots , Q_n\in q_n$ and $P^{\prime}\in p$ such that  $P^{\prime}$ and $Q_1\cup\dots \cup Q_n$ have asymptotically disjoint neighborhoods and 
$Q\setminus (Q_1\cup\dots \cup Q_n)$
 is bounded. 
It follows  that  $P^{\prime}, Q$ have disjoint asymptotic  neighborhoods.

At last,  we find $p_1, \dots ,  p_m\in  \mathcal{B}_{ (X, \mathcal{E})} ^{\sharp}$ and $P_1\in p_1, \dots ,  P_m\in p_m$ such that  $p_1 \cup  \dots  \cup p_m $ and  $Q$  have  disjoint asymptotic   neighborhoods and $P\setminus (P_1 \cup  \dots  \cup P_m )$ is bounded. 
Hence, $P$ and $Q$ have  disjoint asymptotic  neighborhoods. $ \ \Box$

%\bibliography{mybibfile}

\vskip 5pt

CONTACT INFORMATION

I.~Protasov: \\
Faculty of Computer Science and Cybernetics  \\
        Kyiv University  \\
         Academic Glushkov pr. 4d  \\
         03680 Kyiv, Ukraine \\ i.v.protasov@gmail.com

\medskip

K.~Protasova:\\
Faculty of Computer Science and Cybernetics \\
        Kyiv University  \\
         Academic Glushkov pr. 4d  \\
         03680 Kyiv, Ukraine \\ ksuha@freenet.com.ua

\end{document}